\documentclass[11pt]{amsart}
\usepackage{amsmath,amsfonts,amstext,amsthm,epsfig}

\newtheorem{thm}{Theorem}[section]

\newtheorem{lem}[thm]{Lemma}

\theoremstyle{remark}
\newtheorem*{rem}{Remark}

\numberwithin{equation}{section}

\newcommand{\Real}{\mathbb R}

\newcommand{\B}{\mathcal{B}}
\newcommand{\weakto}{\rightharpoonup}

\def\N{\mathbb N}
\def\R{\mathbb R}

\begin{document}

\title[Reaction-diffusion equations in 
unbounded domains]{A remark on reaction-diffusion equations\\ 
in unbounded domains
}%
\author[Martino Prizzi]{}

\email{prizzi@mathsun1.univ.trieste.it}%

\subjclass{35K57, 35B40, 35B41}%
\keywords{Reaction-diffusion equation, asymptotic compactness, attractor}%


\maketitle

\centerline{\scshape  Martino Prizzi}
\medskip

{\footnotesize
\centerline{ Dipartimento di Scienze Matematiche }
\centerline{ Universit\`a degli Studi di Trieste}
\centerline{ Via Valerio 12/b, 34100 Trieste, Italy }
}
\medskip

\bigskip


\begin{quote}{\normalfont\fontsize{8}{10}\selectfont
{\bfseries Abstract.}
We prove the existence of a compact $(L^2-H^1)$ attractor for a 
reaction-diffusion equation in $\Real^N$. This improves a previous 
result of B. Wang
concerning the existence of a compact $(L^2-L^2)$ attractor for the 
same equation. 
\par}
\end{quote}


\section{Introduction} 

In this note we consider the reaction-diffusion equation
\begin{equation}\label{eq}
u_t=\Delta u-\lambda u+f(u)+g(x)\quad\text{in $]0,+\infty[\times\R^N$}.
\end{equation}
We make the following assumptions:
\begin{equation}\label{cond1}
\lambda>0,\quad g\in L^2(\R^N);
\end{equation}
\begin{equation}\label{cond2}
f(0)=0;\quad f(s)s\leq 0,\quad f'(s)\leq C \quad\text{for all $s\in\R$};
\end{equation}
\begin{equation}\label{cond3}
|f'(s)|\leq C(1+|s|^{\beta})\quad\text{for all $s\in\R$},
\end{equation}
where $C$ is some positive constant and
\begin{equation}\label{cond4}
0\leq\beta\quad\text{if $N\leq2$};\quad 0\leq\beta\leq\min\{(2^*/2)-1,4/N\}
\quad\text{if $N\geq3$}.
\end{equation}

In the sequel, we write $L^2$ and $H^1$ instead of 
$L^2(\R^N)$ and $H^1(\R^N)$. Moreover, we denote 
by $\langle\cdot,\cdot\rangle$
the scalar product of $L^2$, by $\|\cdot\|$ the norm of 
$L^2$ and by $\|\cdot\|_E$ the norm of any other space $E$.
It is well known (see e.g. \cite{babinviscik2})
that the Cauchy problem for equation (\ref{eq}) is well posed for initial data
\begin{equation}\label{incon}
u(0,x)=u_0(x) 
\end{equation}
in the space $L^2$.
Actually, equation (\ref{eq}) generates a global semiflow
$\pi$ in $L^2$ and in $H^1$. The following lemma summarises some well 
known a-priori estimates for the semiflow $\pi$ 
(see \cite{wang}, Lemmas 1 and 2):

\begin{lem}\label{priori} Assume conditions (\ref{cond1})-(\ref{cond4}) 
are satisfied. Then

\begin{enumerate}
\item there exists a constant $K>0$ and for every $R>0$ there exists $T(R)>0$ 
such that, if $u\in L^2$ and $\|u\|\leq R$, then $\|\pi(t,u)\|_{H^1}\leq K$ 
for all
$t\geq T(R)$;
\item For every $R>0$ there exists a constant $K(R)$ such that, if $u\in H^1$
and $\|u\|_{H^1}\leq R$, then $\|\pi(t,u)\|_{H^1}\leq K(R)$ for all $t\geq 0$.

\end{enumerate}
\end{lem}

It is a classical result that, if $\R^N$
is replaced by  a bounded domain 
$\Omega\subset\R^N$, equation (\ref{eq}) possesses a compact 
$(L^2-H^1)$ attractor (see e.g. \cite{babinviscik}, \cite{hale} 
and \cite{henry}). The proof of this fact relies essentially on the 
compactness of the Sobolev embedding $H^1(\Omega)\subset L^2(\Omega)$.
If $\Omega$ is unbounded, 
the Sobolev embedding may be no longer compact and new techniques are 
needed. In \cite{babinviscik2} Babin and Vishik  proved the existence of a
compact attractor for equation (\ref{eq}) in an unbounded domain, but 
they needed to introduce weighted spaces in order to overcome the difficulties
arising from the lack of compactness. 
The choice of weighted spaces, however, 
imposes some severe conditions on the forcing term $g$ and on 
the initial data. In \cite{wang}, B. Wang obtained for the first time 
the existence of a compact attractor for equation (\ref{eq}) in the `natural' 
space $L^2$. The crucial step in his proof is the
following 

\begin{thm}[Wang '98]\label{wang1} 
Let $(u_n)_{n\in\N}$ be a bounded sequence in $L^2$ and let $(t_n)_{n\in\N}$
be a sequence of positive numbers, $t_n\to\infty$ as $n\to\infty$.
Then there exists a strictly increasing sequence of natural numbers
$(n_k)_{k\in\N}$ and a function $u\in L^2$ such that $\pi(t_{n_k}, u_{n_k})\to
u$ in $L^2$ as $k\to\infty$. In other words, $\pi$ is asymptotically compact
in the strong $L^2$ topology.
\end{thm}

Combining Lemma \ref{priori} and Theorem \ref{wang1} 
with the abstract results of \cite{hale}, \cite{Lady} and \cite{temam}, 
one easily obtains

\begin{thm}[Wang '98]\label{wang2}
The semiflow $\pi$ possesses a compact $(L^2-L^2)$ attractor.
\end{thm}

On the other hand, the `natural' energy associated to 
equation (\ref{eq}) is given by
\begin{equation}
{\mathcal E}(u)=\frac12\int_{\R^N}|\nabla u|^2\, dx+
\frac\lambda2\int_{\R^N}|u|^2\,dx-\int_{\R^N}F(u)\,dx-\int_{R^N}gu\,dx,
\end{equation}
where $F'(s)=f(s)$, $s\in\R$. The gradient structure of equation (\ref{eq})
implies that, along any given nonconstant solution of (\ref{eq}), the energy 
${\mathcal E}$ decays  and
asymptotically approaches some limit level. The question then arises, whether 
the $\omega$-limit of a given orbit lies at the energy level asymptotically
approached by the orbit itself.
Since ${\mathcal E}$ is continuous on $H^1$, the problem becomes that of proving the
asymptotic compactness of $\pi$ in the strong $H^1$ topology. However, it seems
definitely not trivial to obtain such a result by means of estimates involving
the energy ${\mathcal E}$.

Very recently Efendiev and Zelik in \cite{efzel} 
considered the more general problem
\begin{equation}\label{eqgrad}
u_t=\Delta u-\lambda u+f(u,\nabla u)+g(x)
\end{equation}
in {\em three} spatial dimensions. Using energy estimates and comparison 
arguments
they proved, among other things, the existence of a compact 
$(H^\alpha-H^\alpha)$ attractor
for (\ref{eqgrad}), where $H^\alpha$ is some fractional space between 
$H^2$ and $H^1$. Their technique, however, exploits the Sobolev imbedding
$H^2\subset L^\infty$, which is no longer true in higher space dimensions. 

We have found more  
convenient to follow a different way: in this note we show that the asymptotic
compactness in $H^1$ can be recovered from the asymptotic compactness
in $L^2$ by a simple continuity argument (Theorem \ref{ascomp}). 
As a consequence, 
we deduce that the $(L^2-L^2)$ attractor, whose existence was 
established by B. Wang, is actually an $(L^2-H^1)$ attractor, like in the case 
of bounded domains (Theorem \ref{attratt}). 
The proof is very simple and is based on Henry's theory
of abstract parabolic equations (see \cite{henry}). 

Finally, we would like to mention the recent results obtained by 
Zelik in \cite{zel1}, concerning the existence and the entropy of
{\em locally compact} attractors for equations like (\ref{eq}). 
See also \cite{zel2}, where similar results have 
been obtained for a dampded wave equation in $\R^N$.

\section{The proof}

Define the bilinear form
\begin{equation}
a(u,v):=\int_{\R^N}\nabla u(x)\cdot\nabla v(x) dx,\quad u,v\in H^1.
\end{equation}
The form $a(\cdot,\cdot)$ uniquely determines a self-adjoint operator
$A\colon D(A)\subset L^2\to L^2$, defined by the relations
\begin{equation}
\begin{cases}
D(A):=\{u\in H^1\mid\text{$\exists\, w\in L^2$ such that
$a(u,v)=\langle w,v\rangle$ for all $v\in H^1$}\}\\
Au:=w\quad\text{for $u\in D(A)$}
\end{cases}
\end{equation}
The self-adjoint operator $-A$ is the generator of an analytic semigroup 
of linear operators $e^{-At}$, $t\geq 0$, satisfying the following estimates:
\begin{equation}
\|e^{-At}u\|\leq M e^{at}\|u\|\quad\text{for all $u\in L^2$ and all $t\geq 0$},
\end{equation}
\begin{equation}
\|e^{-At}u\|_{H^1}\leq 
M e^{at}t^{-1/2}\|u\|\quad\text{for all $u\in L^2$ and all $t> 0$},
\end{equation}
where $M$ and $a$ are two positive constants. We need also the following
\begin{lem}\label{nemitski}
The assignement $u\mapsto f\circ u$ defines a map $\hat f\colon H^1\to L^2$,
which is Lipschitz continuous on every bounded set in $H^1$.
\end{lem}
\begin{proof}
By (\ref{cond3}) we have
$$
|f(s)|\leq C(|s|+|s|^{\beta+1}).
$$
Let $u\in H^1$. Then 
\begin{multline*}
\int_{\R^N}|f(u(x)|^2dx 
\leq 2C\left(\int_{\R^N}|u(x)|^2dx+
\int_{\R^N}|u(x)|^{2(\beta+1)}dx\right)\\
=2C\left(\|u\|^2+\|u\|^{2(\beta+1)}_{L^{2(\beta+1)}}\right).\\
\end{multline*}
Since $2\leq 2(\beta+1)\leq 2^*$, $H^1$ is continuously imbedded in 
$L^{2(\beta+1)}$. Then there exists a constant $C_1$ such that
$$
\int_{\R^N}|f(u(x)|^2dx 
\leq 2C\left(\|u\|^2+C_1\|u\|_{H^1}^{2(\beta+1)}\right)
$$
and hence
$$
\|\hat f(u)\|\leq 2C\left(\|u\|+C_1\|u\|_{H^1}^{(\beta+1)}\right).
$$
This shows that $\hat f\colon H^1\to L^2$ is well defined and maps
bounded sets of $H^1$ to bounded sets of $L^2$.

Now let $s_1,s_2\in\R$. By the mean value theorem and by (\ref{cond3}) we have
\begin{multline*}
|f(s_1)-f(s_2)|=|f'(\theta s_1+(1-\theta)s_2)(s_1-s_2)|\\
\leq C(1+|\theta s_1+(1-\theta)s_2|^\beta)|s_1-s_2|
\leq C_1(1+|s_1|^\beta+|s_2|^\beta)|s_1-s_2|,
\end{multline*}
where $\theta=\theta(s_1,s_2)$, $0\leq\theta\leq1$, 
and $C_1$ is a positive constant. 
Let $u_1,u_2\in H^1$, $\|u_1\|_{H^1},\|u_2\|_{H^1}\leq R$. Then we have
\begin{multline*}
\int_{\R^N}|f(u_1(x))-f(u_2(x))|^2dx\\
\leq 2C_1\int_{\R^N}(1+
|u_1(x)|^{2\beta}+|u_2(x)|^{2\beta})|u_1(x)-u_2(x)|^2dx.\\
\end{multline*}
Let $p:=2^*/(2\beta)>1$, $q:=p/(p-1)=2^*/(2^*-2\beta)$. Then, by H\"older 
inequality, we get
\begin{multline*}
\int_{\R^N}|f(u_1(x))-f(u_2(x))|^2dx
\leq 2C_1\int_{\R^N}|u_1(x)-u_2(x)|^2dx\\
+ 2C_1\left(\|u_1\|_{L^{2^*}}^{2\beta}+\|u_2\|_{L^{2^*}}^{2\beta}\right)
\|u_1-u_2\|_{L^{2q}}^2.
\end{multline*}
Since $2\leq 2q\leq 2^*$, $H^1$ is continuously imbedded in $L^{2q}$. Then 
there exists a constant $C_2$ such that
$$
\int_{\R^N}|f(u_1(x))-f(u_2(x))|^2dx\leq 
C_2(1+\|u_1\|_{H^1}^{2\beta}+\|u_2\|_{H^1}^{2\beta})\|u_1-u_2\|_{H^1}^2
$$
and hence
$$
\|\hat f(u_1)-\hat f(u_2)\|\leq C_2^{1/2}(1+2R^{2\beta})^{1/2}
\|u_1-u_2\|_{H^1}.
$$
The proof is complete.
\end{proof}

Following Henry \cite{henry}, the Cauchy problem (\ref{eq})-(\ref{incon}), 
for initial data  $u_0\in H^1$, can be formulated as an
abstract parabolic initial value problem
\begin{equation}\label{abstreq}
\begin{cases}
\dot u+Au=-\lambda u+\hat f(u)+g\\
u(0)=u_0
\end{cases}
\end{equation}

It is well known that equation (\ref{abstreq}) is 
equivalent to the integral equation
\begin{equation}\label{inteq}
u(t)=e^{-At}u_0+\int_0^t e^{-A(t-s)} (-\lambda u +\hat f(u(s))+g)ds.
\end{equation}
We have the following crucial

\begin{lem}\label{singron}
Let $(u_{n})_{n\in\N}$ be a sequence in $H^1$, let $u\in H^1$, and assume
that $u_n\weakto u$ in $H^1$, $u_n\to u$ in $L^2$. Then $\pi(t,u_n)\to\pi(t,u)$
in $H^1$, uniformly on $[t_0,t_1]$, for all $t_1>t_0>0$.
\end{lem}
\begin{proof}
Let $t_1>0$ be fixed. Since the set $\{u_n\mid n\in\N\}\cup\{u\}$ is bounded
in $H^1$, by Lemma \ref{priori} there exists $R>0$ such 
that $\|\pi(t,u_n)\|_{H^1}\leq R$
and $\|\pi(t,u)\|_{H_1}\leq R$, for all $t\in[0,t_1]$ and for 
all $n\in\N$. Let $L$ be a 
Lipschitz constant for $\hat f$ on the ball of radius $R$ in $H^1$. 
Write $u_n(t):=\pi(t,u_n)$ and $u(t):=\pi(t,u)$.
For $t\in[0,t_1]$, we have
$$
u_n(t)=e^{-At}u_n+\int_0^t e^{-A(t-s)} (-\lambda u_n(s)+\hat f(u_n(s))+g)ds
$$
and
$$
u(t)=e^{-At}u+\int_0^t e^{-A(t-s)} (-\lambda u +\hat f(u(s))+g)ds.
$$
It follows that, for $t\in[0,t_1]$,
\begin{multline*}
\|u_n(t)-u(t)\|_{H^1}\leq Me^{at_1}t^{-1/2}\|u_n-u\|\\+M(\lambda+L)e^{at_1}
\int_0^t(t-s)^{-1/2}\|u_n(s)-u(s)\|_{H^1}ds.\\
\end{multline*}
By the singular Gronwall's Lemma (see \cite{henry}), there is a constant
$C_1=C_1(\lambda,L,M,a,t_1)$ such that, for $t\in[0,t_1]$,
\begin{multline*}
\|u_n(t)-u(t)\|_{H^1}\leq Me^{at_1}t^{-1/2}\|u_n-u\|\\
+C_1\int_0^t(1+(t-s)^{-1/2})Me^{at_1}s^{-1/2}\|u_n-u\|ds.\\
\end{multline*}
Finally, a simple computation shows that there exists a constant 
$C_2=C_2(\lambda,L,M,a,t_1)$ such that, for $t\in[0,t_1]$,
$$
\|u_n(t)-u(t)\|_{H^1}\leq C_2(1+t^{-1/2})\|u_n-u\|.
$$
This completes the proof.
\end{proof}

\begin{rem}
It seems that such a result cannot be obtained as a consequence of
a-priori estimates like the ones in \cite{wang}. In a different context, a
similar observation was made also in \cite{pr}.
\end{rem}

The next theorem shows that, thanks to Lemma \ref{singron}, the asymptotic 
compactness of $\pi$ in $L^2$ implies the asymptotic compactness of $\pi$ in 
$H^1$.

\begin{thm}\label{ascomp}
Let $(u_n)_{n\in\N}$ be a bounded sequence in $H^1$ and let $(t_n)_{n\in\N}$
be a sequence of positive numbers, $t_n\to\infty$ as $n\to\infty$.
Then there exists a strictly increasing sequence of natural numbers
$(n_k)_{k\in\N}$ and a function $u\in H^1$ such that $\pi(t_{n_k}, u_{n_k})\to
u$ in $H^1$ as $k\to\infty$. In other words, $\pi$ is asymptotically compact
in the strong $H^1$ topology.
\end{thm}
\begin{proof}
Fix any positive $T$. Since $t_n\to\infty$ as $n\to\infty$, we have $t_n>T$
for all sufficiently large $n$, say $n\geq n_0$.
Since the sequence $(u_n)_{n\in\N}$ is bounded in $H^1$, by Lemma \ref{priori}
also the sequence $\pi(t_n-T,u_n)_{n\geq n_0}$ is bounded in $H^1$. Then 
there exists a strictly increasing sequence of natural numbers
$(n_k)_{k\in\N}$, $n_k\geq n_0$ for all $k\in\N$, 
and a function $\bar u\in H^1$ such that $\pi(t_{n_k}-T, u_{n_k})
\weakto \bar u$ in $H^1$ as $k\to\infty$. On the other hand, by Theorem 
\ref{wang1}, we 
can choose the sequence $(n_k)_{k\in\N}$ in such a way that 
$\pi(t_{n_k}-T, u_{n_k})\to \bar u$ in $L^2$ as $k\to\infty$. 
Then, by Lemma \ref{singron}, we have
$$
\pi(t_{n_k},u_{n_k})=\pi(T,\pi(t_{n_k}-T,u_{n_k}))\to\pi(T,\bar u)\quad
\text{in $H^1$ as $k\to\infty$.}
$$
The proof is complete.
\end{proof}

Finally, we can state and prove

\begin{thm}\label{attratt}
The semiflow $\pi$ possesses a compact $(L^2-H^1)$ attractor.
\end{thm}

\begin{proof}
By Lemma \ref{priori}  there exists a bounded set $\B$ in $H^1$ and for
any bounded set $B$ in $L^2$ there exists $T(B)>0$ 
 such that $\pi(t,B)\subset\B$ for all $t\geq T(B)$. On the other hand, 
by Theorem \ref{ascomp} $\pi$ is asymptotically compact in $H^1$. The 
conclusion follows 
from the abstract results of \cite{hale}, \cite{Lady} and \cite{temam}.
\end{proof}

\begin{rem}Of course, by
Lemma \ref{priori}, the $(L^2-L^2)$ attractor and the $(L^2-H^1)$ 
attractor coincide.
Moreover, by the continuity result of Lemma \ref{singron}, the 
Hausdorff (or the fractal) dimension of the attractor is the same in
$L^2$ and in $H^1$.
\end{rem}

\smallskip


\end{document}